\documentclass[a4paper, 11pt]{amsart}
\usepackage{amssymb, amsmath, graphicx}
\newtheorem{theorem}{Theorem}
\newtheorem*{theorem*}{Theorem}
\newtheorem{lemma}[theorem]{Lemma}

\newtheorem{proposition}[theorem]{Proposition}
\newtheorem{conjecture}{Conjecture}
\newtheorem*{conjecture*}{Conjecture}

\theoremstyle{definition}

\theoremstyle{remark}

\newcommand{\diam}{\textup{diam}}

\begin{document}

\title[A note on the Hayman-Wu theorem]{A note on the Hayman-Wu theorem}
\date{\today}

\author[E. Crane]{Edward Crane}
\email{edward.crane@gmail.com}
\address{Mathematics Department\\
         University of Bristol\\
         University Walk\\ Bristol  BS8 1TW\\
         United Kingdom}


\begin{abstract}
The Hayman-Wu theorem states that the preimage of a line or circle $L$ under a conformal mapping from the unit disc $\mathbb{D}$ to a simply-connected domain $\Omega$ has total Euclidean length bounded by an absolute constant. The best possible constant is known to lie in the interval $[\pi^2, 4\pi)$, thanks to work of {\O}yma and Rohde. Earlier, Brown Flinn showed that the total length is at most $\pi^2$ in the special case in which $L \subset \Omega$. Let $r$ be the anti-M\"{o}bius map that fixes $L$ pointwise. In this note we extend the sharp bound $\pi^2$ to the case where each connected component of $\Omega \cap r(\Omega)$ is bounded by one arc of $\partial\Omega$ and one arc of $r(\partial \Omega)$. We also strengthen the bounds slightly by replacing Euclidean length with the strictly larger spherical length on $\mathbb{D}$.
\end{abstract}

\maketitle

\section{Introduction}
Let $g$ be a conformal mapping from the unit disc $\mathbb{D}$ onto a simply-connected domain $\Omega$ in the Riemann sphere. Let $L$ be any line or circle. Let $\Lambda_1$ denote one-dimensional Hausdorff measure with respect to the Euclidean metric in $\mathbb{C}$. The Hayman-Wu theorem \cite{HW} says that there exists a constant $C$, which does not depend on $g$ and $L$, such that $\Lambda_1\left(g^{-1}(L)\right) \le C$. That is, the sum of the Euclidean lengths of all components of $g^{-1}(L)$ is bounded by an absolute constant. Following Rohde we denote the smallest such constant by $\text{\O}$ in memory of Knut {\O}yma who gave in \cite{Oy1} a simple proof that $\text{\O} \le 4 \pi$ and then in \cite{Oy2} constructed a family of examples to show that $\text{\O} \ge \pi^2$.  A concise exposition of these proofs together with the required background may be found in \cite{GM}. Rohde later showed in \cite{R} that $\text{\O} < 4 \pi$ by combining {\O}yma's method with an inequality relating the geodesic curvature of $L$ to the change of hyperbolic density when passing from $f(\mathbb{D})$ to a certain subdomain. He did not give an explicit bound less than $4\pi$, but this could in principle be done using his method.

 Earlier, in \cite{BF}, Barbara Brown Flinn showed that the upper bound $\pi^2$ is valid in the special case in which $L \subset \Omega$; she gave a bound on the Euclidean length of the perimeter of a hyperbolically convex subset of $\mathbb{D}$ and then applied a result of Vilhelm J{\o}rgensen \cite{J} to the effect that a disc is a hyperbolically convex subset of any hyperbolic domain that contains it. This shows that the Jordan curve $g^{-1}(L)$ bounds a hyperbolically convex subset of $\mathbb{D}$ and therefore has length less than $\pi^2$.

In this note we extend the above result to give the sharp bound $\pi^2$ for a more general class of pairs $(\Omega, L)$. This is precisely the class in which only the simplest special case arises in {\O}yma's proof of the bound $\textup{\O} \le 4\pi$, although in fact our proof is apparently unrelated to that of {\O}yma and much closer in spirit to that of Brown Flinn. We also strengthen the result slightly by using the one-dimensional Hausdorff measure with respect to the spherical metric $2|dz|/(1+|z|^2)$, which we denote by $\Lambda_1(\cdot, \sigma)$.

\begin{theorem}\label{T: main}{\quad}\\
 Let $g: \mathbb{D} \to \Omega \subset \mathbb{C}_\infty$ be a conformal homeomorphism and let $L$ be a line or circle. Let $r$ be the unique anti-M\"{o}bius map that fixes $L$ pointwise. Suppose that each connected component of $\Omega \cap r(\Omega)$ is bounded by one connected component of $\overline{r(\Omega)} \cap \partial \Omega$ together with its image under $r$. Then $g^{-1}(L)$ bounds a hyperbolically convex subset of $\mathbb{D}$ and hence \[\Lambda_1\left(g^{-1}(L)\right)\, <\, \Lambda_1\left(g^{-1}(L),\sigma\right)\, <\, \pi^2\,.\]
\end{theorem}

By the same method we solve a closely related extremal problem which may be of independent interest. Let $\mathbb{H}$ denote the upper half-plane. If $V$ is a domain in $\mathbb{C}$ bounded by a Jordan curve, and $\gamma$ is an arc of that boundary, we say that $f$ is a \emph{conformal reflection} across $\gamma$ if $f: V \to \mathbb{C} \setminus \overline{V}$ is an analytic injection and the following defines a conformal automorphism $\tilde{f}$ of the domain $V \cup \gamma \cup f(V)$:
\[\tilde{f}(z) \,:=\, \begin{cases}f(z) & \text{if $z \in V$}, \\ z & \text{if $z \in \gamma$},\\ f^{-1}(z) & \text{if $z \in f(V)$}.   \end{cases}\] 

\begin{proposition}\label{P: conformal reflection}{\quad}\\
Let $\gamma$ be a simple curve in $\mathbb{H}$ that lands at $-1$ and $1$. Then $\gamma$ together with the interval $[-1,1]$ bounds a simply-connected domain $B \subset \mathbb{H}$. Suppose that there exists a conformal reflection $f$ across $\gamma$ that maps $B$ inside $\mathbb{H}$. Then $\gamma$ is a real-analytic curve of Euclidean length at most $\pi$, with equality if and only if $\gamma$ is a semicircle on the diameter $[-1,1]$.
\end{proposition}

We now give two alternative statements of essentially the same result.

\begin{proposition}\label{P: totally real schlicht functions}{\quad}\\
 Let $g: \mathbb{D} \to \mathbb{C}$ be a totally real univalent function, i.e. in the Taylor expansion $g(z) = \sum a_n z^n$, all the $a_n$ are real, and suppose that $a_1 > 0$. Suppose that $g$ has finite radial limits at $-1$ and $1$, say $\lim_{z \to -1}g(z) = a$ and $\lim_{z \to 1} g(z) = b$. Let $\Gamma$ be the circular arc of points $w \in \mathbb{D}$ such that $\arg\left(w - 1)/(w+1)\right) = 3\pi/4$. Then $g(\Gamma)$ is a curve of Euclidean length at most $\pi(b-a)$, with equality if and only if $g(\mathbb{D})$ has the form $\mathbb{C} \setminus ((-\infty,a] \cup [b,\infty))$. 
\end{proposition}

In the following, $\omega(z, E, U)$ denotes the harmonic measure of a set $E \subset \partial U$ with respect to a point $z$ in the domain $U$.

\begin{proposition}\label{P: level set of harmonic measure}{\quad}\\
Let $U$ be a subdomain of $\mathbb{H}$ and denote by $\partial U$ the boundary of $U$ in the Riemann sphere $\mathbb{C}_\infty$. Suppose that $[-1,1]
 \subset \partial U$ and $\partial U \setminus (-1,1)$ is connected. For $0 < \alpha < 1$ we define the level set 
 \[ \gamma_\alpha \, = \, \{ z \in U \, : \, \omega(z, (-1,1), U) = \alpha \}\,. \]
Then $\Lambda_1\left(\gamma_{1/2}\right) \, \leqslant \, \pi$, with equality if and only if $U = \mathbb{H}$.
\end{proposition}

\noindent{Remark:} The method of {\O}yma \cite{Oy1} can be applied to the situations of these propositions, in each case giving a slightly weaker length bound with $4$ in place of $\pi$. If we had been able to find a counterexample to Proposition~\ref{P: conformal reflection}, then we would have been able to make a construction along the lines of \cite{Oy2} to show that $\text{\O} > \pi^2$.

The formulation of Proposition~\ref{P: level set of harmonic measure} naturally leads to the following conjecture.

\begin{conjecture}
Let $U$ and $\gamma_{\alpha}$ be defined as in Proposition~\ref{P: level set of harmonic measure}. Then 
\[ \Lambda_1(\gamma_\alpha) \, \leqslant \, \frac{2 \pi (1-\alpha)}{\sin \pi \alpha} \,,  \]
with equality if and only if $U = \mathbb{H}$.
\end{conjecture}

\section{Hyperbolic convexity and spherical length}
 Let $U$ be a simply-connected hyperbolic domain, and denote by $\rho_U$ its associated complete hyperbolic metric. A subset $E \subset U$ is said to be \emph{hyperbolically convex} if for each $z,w \in E$, the unique $\rho_U$-geodesic segment joining $z$ to $w$ is contained in $E$. In the special case where $U$ is the unit disc, we can apply the transformation $K: z \mapsto 2z/(1+|z|^2)$, which is an isometry from the Poincar\`{e} model of the hyperbolic plane to the Klein model of the hyperbolic plane. In the Klein model, the geodesics are chords of the unit circle. Thus the hyperbolically convex subsets of $\mathbb{D}$ are precisely those which map under $K$ to ordinary convex subsets of plane. It is well-known that any bounded convex subset $E'$ of the Euclidean plane has a rectifiable boundary, whose length satisfies the isoperimetric inequalities 
\[2 \,\diam_{\textup{Euc}}(E')\, \le\, \Lambda_1(\partial E')\, \le\, \pi\, \diam_{\textup{Euc}}(E')\,.\]
Here $\Lambda_1$ denotes $1$-dimensional Hausdorff measure with respect to the Euclidean metric on $\mathbb{C}$ and $\textup{diam}_{\textup{Euc}}$ is the diameter in the Euclidean metric. The upper bound follows from Santal\`{o}'s integral formula for $\Lambda_1(\partial E)$. Setting $E' = K(E)$ it follows that any hyperbolically convex subset $E$ of $\mathbb{D}$ also has rectifiable boundary. In the following, $[z,w]_\mathbb{D}$ denotes the geodesic segment from $z$ to $w$ in the hyperbolic metric on $\mathbb{D}$. 

\begin{theorem}[Brown Flinn]{\quad}\label{T: Brown Flinn estimate}
\begin{enumerate}
\item \cite[Lemma 4]{BF} Let $z_1, z_2, \dots, z_m$ be the vertices of a convex hyperbolic polygon in $\mathbb{D}$, taken in order. Suppose that $z_2, \dots, z_{m-1}$ are all contained in the interior of the Euclidean convex hull of $[z_1, z_m]_\mathbb{D}$. Then
\[ \sum_{i=1}^{m-1} \Lambda_1\left( [z_i, z_{i+1}]_\mathbb{D}\right) \, \leqslant \, \Lambda_1\left([z_1, z_m]_\mathbb{D}\right)\,. \] 
\item \cite[Theorem 3]{BF} Let $E$ be a hyperbolically convex subset of $\mathbb{D}$. Then
 \[ \Lambda_1(\partial E) \, \leqslant \, \frac{\pi^2}{2} \, \textup{diam}_{\textup{Euc}}(E) \, \leqslant \, \pi^2\,.\]
\end{enumerate}
\end{theorem}

 Brown Flinn proved part (1) of Theorem~\ref{T: Brown Flinn estimate} as the main tool in the proof of part (2), by an approximation argument. We will use a similar approximation argument to prove an apparently weaker Euclidean length estimate in $\mathbb{H}$ instead of $\mathbb{D}$, which we then use to provide an estimate of the spherical length of $\partial E$.
  
\begin{lemma}\label{L: Brown Flinn estimate in the upper half-plane} {\quad}\\
Let $V$ be a hyperbolically convex subset of $\mathbb{H}$, and let $z, w \in \partial V$, $z \neq w$. Let $\kappa$ be the $\rho_\mathbb{H}$-geodesic that passes through $z$ and $w$. Suppose that $\kappa$ lands at $a, b \in \mathbb{R}$, so that  $\kappa$ and $[a,b]$ together bound an open semidisc $D$, and suppose that $V \cap D$ is contained in the Euclidean convex hull of $[z,w]_\mathbb{H}$. Then \[\Lambda_1(\partial V \cap (D \cup [z,w]_\mathbb{H}))\,\le\, \Lambda_1([z,w]_\mathbb{H}) \,\le\, \frac{\pi}{2} |z - w|\,,\] with equality if and only if $V \cap D = \emptyset$.
\end{lemma} 
\begin{proof}
We check that part (i) of Theorem~\ref{T: Brown Flinn estimate} applies with $\mathbb{D}$ replaced by $\mathbb{H}$. We can think of this as a limiting case when the polygon is close to the unit circle and has small Euclidean diameter. To make this precise, we consider the image $(z_1, \dots, z_m)$ in $\mathbb{D}$ of a hyperbolic polygon $(w_1, \dots, w_m)$ in $\mathbb{H}$ under the M\"{o}bius map $ M_n: w \mapsto (w - ni)/(w+ni)$. This map is an isometry from $\rho_\mathbb{H}$ to $\rho_\mathbb{D}$, and we have $n M_n'(w) \, \to \, -2i$ locally uniformly in $w$, as $n \to \infty$. Suppose that $w_2, \dots, w_{m-1}$ are contained in the interior of the Euclidean convex hull of $[w_1, w_m]_\mathbb{H}$. Then for $n$ sufficiently large $z_2, \dots, z_{m-1}$ are contained in the interior of the Euclidean convex hull of $[z_1, z_m]_\mathbb{D}$. In the limit as $n \to \infty$ we obtain the required inequality. Since the length depends continuously on the $w_i$, the inequality still holds if the $w_i$ are allowed to be in the boundary of the Euclidean convex hull of $[w_1, w_m]_\mathbb{H}$. To obtain the full result we need an approximation argument. Note that for any $\epsilon > 0$ we may choose a convex hyperbolic polygon $w_1, w_2, \dots, w_m$ with successive vertices no more than $\epsilon$ apart in the Euclidean metric, such that $w_1 = z, w_m = w$ and each $w_i \in \partial V$, and such that $\partial V \cap D$ is covered by the union of the closed discs with Euclidean centre $(w_i + w_{i+1})/2$ and Euclidean radius $(1+\epsilon)\,.\,|w_{i+1} - w_i|$, where $i$ runs over $\{1, \dots, m-1\}$. This yields our estimate on the Hausdorff measure $\Lambda_1(\partial V \cap D)$. For equality to occur, it must not be possible to subdivide $\partial(V \cap D)$ into two subarcs meeting at a point $\zeta \in D$, for then we would have \[\Lambda_1(\partial V \cap D)\, \leqslant \, \Lambda_1([z,\zeta]_\mathbb{H})\, +\, \Lambda_1([\zeta,w]_\mathbb{H}) \,<\, \Lambda_1([z,w]_\mathbb{H})\,.\] 
\end{proof}

One cannot reverse the argument with the M\"{o}bius maps $M_n$ to recover Theorem~\ref{T: Brown Flinn estimate} directly from Lemma~\ref{L: Brown Flinn estimate in the upper half-plane}, so this does not appear to constitute progress! However, Lemma~\ref{L: Brown Flinn estimate in the upper half-plane} is exactly what we will need for proving Propositions~\ref{P: conformal reflection},~\ref{P: totally real schlicht functions} and~\ref{P: level set of harmonic measure}, and is also what we need next to estimate the spherical length of the boundary of a hyperbolically convex subset of $\mathbb{D}$. 
 
  The spherical metric on $\mathbb{C}_\infty$ is the metric $\sigma$ associated to the conformal metric $2\,|dz|/(1+|z|^2)$. We denote the one-dimensional Hausdorff measure of $X \subset \mathbb{C}_\infty$ with respect to the metric $\sigma$ by $\Lambda_1(X, \sigma)$. By comparing the length elements $|dz|$, $2 |dz| / (1+|z|^2)$, and $2 |dz|/(1 - |z|^2)$, we see that for any $z, w \in \mathbb{D}$ we have 
 \[  |z-w| \, \leqslant \, \sigma(z,w) \, \leqslant \, \rho_\mathbb{D}(z,w)\,,\]
 with equality only if $z=w$.
 
\begin{lemma}\label{L: spherical length estimate}{\quad}\\
Let $E$ be a hyperbolically convex subset of $\mathbb{D}$, with boundary $\partial E$ relative to $\mathbb{C}$. Then $\Lambda_1(\partial E, \sigma)\, \leqslant \, \pi^2$, with equality if and only if $E = \mathbb{D}$.
\end{lemma} 
\begin{proof}
To prove the inequality it suffices to deal with the case in which $E$ is bounded by a hyperbolically convex polygon $(z_1, \dots, z_m)$ with non-empty interior and all vertices in $\mathbb{D}$; the approximation argument is as in the proof of Lemma~\ref{L: Brown Flinn estimate in the upper half-plane} and we will not repeat it. 

 Next, we reduce to the case in which $\overline{E}$ contains $0$. If this is not the case, then let $r = \rho_\mathbb{D}(0, E)$; then there is a unique point $z_0$ of $\partial{E}$ such that $\rho_\mathbb{D}(0,z_0) = r$, because $E$ is hyperbolically convex. Consider the M\"{o}bius map \[M: z \mapsto \frac{z - z_0}{1 - z_0^*z}\,.\] Here $z_0^*$ is the complex conjugate of $z_0$. Because $M$ is a hyperbolic isometry, $M(E)$ is hyperbolically convex. $M(z_0) = 0$ so $0 \in \partial(M(E))$. We claim that \[\Lambda_1(\partial M(E), \sigma)\, > \,\Lambda_1(\partial E, \sigma)\,,\] because $\partial E$ is contained in the region where the spherical derivative of $M$ is greater than $1$. This region is bounded by the arc of points $w$ such that the ratio of spherical and hyperbolic densities takes the same value at $w$ and $M(w)$. Since this ratio only depends on $|w|$, this isometric circle is also the isometric circle for $M$ with respect to the Euclidean metric. Explicitly, it is a hyperbolic geodesic which is the perpendicular bisector of $z_0$ and $0$, and it is an arc of the Ford circle of $M$. By replacing $E$ by $M(E)$ if necessary, we may assume that $0 \in \overline{E}$.
 
 We now consider the stereographic projection $\Pi: \mathbb{C} \to S^2$ from the point $(0,0,-1)$. Here $S^2$ is the unit sphere in $\mathbb{R}^3 = \mathbb{C} \times \mathbb{R}$. Explicitly, 
 \[\Pi(z) \,=\, \left(\frac{2z}{1+|z|^2}\,,\, \frac{1 - |z|^2}{1 + |z|^2}\right)\,.\]
This maps $\mathbb{D}$ conformally onto the upper hemisphere of $S^2$. The spherical length of $\partial E$ is the Euclidean length of the image $\Pi(\partial E)$ on $S^2$.  Note that $\Pi$ followed by orthogonal projection back onto $\mathbb{C}$ gives the transformation $K$ described earlier. Define $w_i = \Pi(z_i)$ for $i = 1, \dots , m$, with co-ordinates $(x_i, t_i) \in \mathbb{C} \times \mathbb{R}$.  Note that $K\left([z_i, z_{j}]_\mathbb{D}\right)$ is the Euclidean line segment $[x_i, x_j]$, so $K(\partial E)$ is a Euclidean convex polygon.  This situation is illustrated in figure~\ref{F: stereographic projection} in the case where $L$ is a hyperbolic triangle.

\begin{figure}
	\begin{center}
		\includegraphics[scale=1.80]{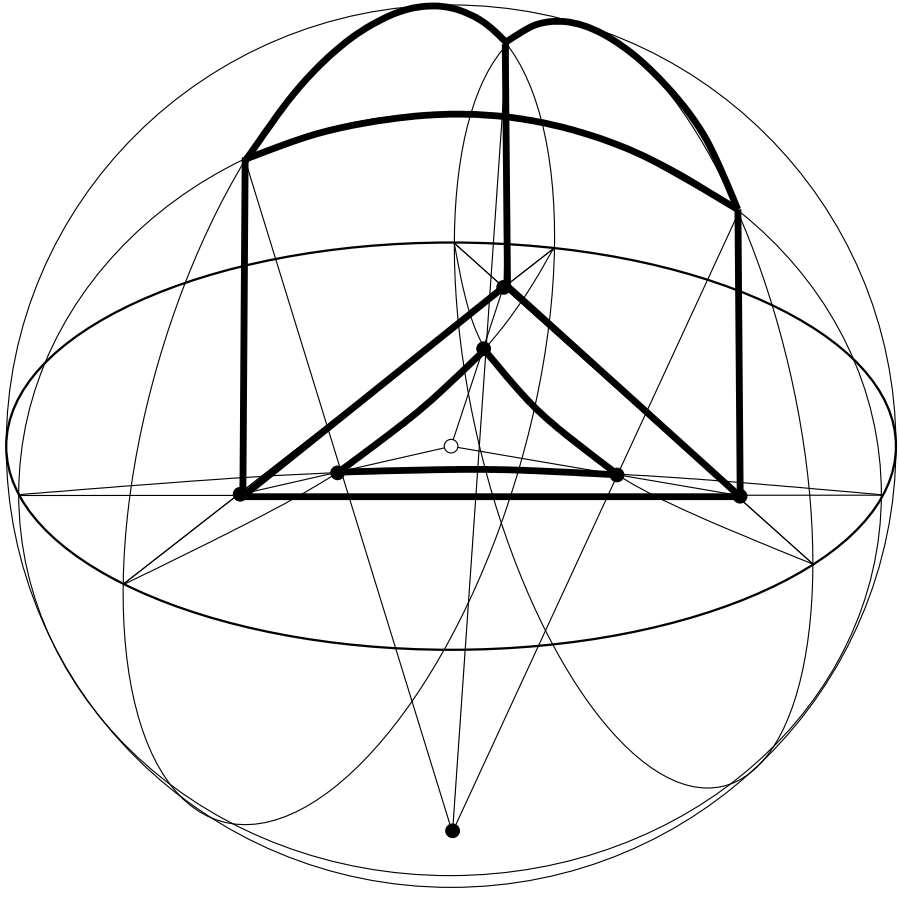}
	\end{center}
\caption{Stereographic projection of a hyperbolic polygon}
	\label{F: stereographic projection}
\end{figure}

Without loss of generality, we may assume that the vertices have been labelled cyclically so that $t_1 = \min(t_1, \dots, t_m)$. Now let $L$ denote the Euclidean length of $K(\partial E)$. We have $L < 2\pi$. Let $\gamma: [0,L] \to K(\partial E)$ be a parametrization by Euclidean arc length, so that $\gamma(0) = K(z_1) = \gamma(L)$.  We can then define a mapping \[\Phi: [0,L] \times (0,\infty) \,\to\, K(\partial E) \times (0,\infty)\] by \[ \Phi(t + is) = (\gamma(t),s)\,.\] 
This mapping \emph{unrolls} the piecewise flat surface $K(\partial E) \times (0, \infty)$ and is a local isometry for the Euclidean metric on each flat piece. Thus we have 
\[\Lambda_1( \partial E, \sigma) \, = \, \Lambda_1(\Pi(\partial E)) \, = \, \Lambda_1\left(\Phi^{-1}(\Pi(\partial E))\right)\,.\]    
The curve $\Gamma = \Phi^{-1}(\Pi(\partial E))$ is a piecewise geodesic curve in $\mathbb{H}$. To see this, note that $\Pi([z_i, z_{i+1}]_\mathbb{D})$ is an arc of a semicircle in $\mathbb{R}^3$ that meets the plane $\mathbb{C} \times \{0\}$ orthogonally at two points on the unit circle.  Let $V$ be the region of $\mathbb{H}$ bounded by $\Gamma$ and the two vertical geodesics $[i t_1, \infty]_\mathbb{H}$ and $[L + i t_1, \infty]_\mathbb{H}$. We will show that $V$ is hyperbolically convex. We can then apply Lemma~\ref{L: Brown Flinn estimate in the upper half-plane} with $z = it_1$ and $w = L+ it_1$ to deduce the desired inequality \[\Lambda_1(\partial E, \sigma) \, \leqslant \, \frac{\pi}{2}\, \left|\left(L+it_1\right) - it_1\right| \,  =\, \pi L/2  \,<\, \pi^2\,. \]  

To show that $V$ is convex, we consider two consecutive edges $e_1$ and $e_2$ of $\partial K(E)$. If they meet at $0$, then the lifts $\Phi^{-1}(\Pi(e_1))$ and $\Phi^{-1}(\Pi(e_2))$ are abutting pieces of the same geodesic in $\mathbb{H}$. Otherwise, a calculation shows that away from the discontinuities where $\gamma(t)$ is a vertex of $K(\partial E)$ we have
\[\sin \arg\left( \frac{d}{dt} \Phi^{-1}(\Pi(\gamma(t))) \right) \, = \, |\gamma(t)| \, \cos \arg \left(\frac{\frac{d}{dt}\gamma(t)}{\gamma(t)}\right)\,.\] The jumps of the right-hand side at the discontinuities are all positive because $0 \in \overline{E}$. It follows that at each vertex of $V$ the interior angle is less than or equal to $\pi$, and hence $V$ is hyperbolically convex. 

For the equality statement, we can no longer deal only with the case of polygonal boundary. Instead, we need $\partial K(E)$ to be approximated by a sequence of Euclidean-convex curves in $\mathbb{D}$ whose lengths approach $2\pi$, and therefore approach the unit circle, and it follows that $E = \mathbb{D}$. 
\end{proof}

\section{Proof of Theorem~\ref{T: main}}
 
In the case where $\Omega \,\cap\, r(\Omega)$ is topologically an annulus, which is the case where $L \subset \Omega$, the hyperbolic convexity of the domain bounded by $g^{-1}(L)$ was demonstrated by Brown Flinn in \cite{BF}, and we can apply Lemma~\ref{L: spherical length estimate} to obtain the spherical length estimate.

 Now suppose that $\Omega \cap r(\Omega)$ is not a topological annulus. Then it is the union of countably many connected components $V_k$, each of which is simply-connected. Each component $V_k$ contains a single connected component $L_k$ of $L$. The open sets $U_k  =  g^{-1}(V_k)$ are simply-connected subsets of $\mathbb{D}$, each containing a single arc $\gamma_k = g^{-1}\left(L_k\right)$ of $g^{-1}(L)$. According to the hypothesis, removing $\gamma_k$ separates $U_k$ into two components $A_k$ and $B_k$, where $B_k$ is the component bounded by $\gamma_k$ and an arc of $\partial\mathbb{D}$. We claim that $\mathbb{D} \setminus \cup_{k} B_k$ is a hyperbolically convex subset of $\mathbb{D}$. Hence, by Lemma~\ref{L: spherical length estimate}, the total spherical length of its boundary is less than $\pi^2$.  It suffices to prove that for each $k$ the set $\mathbb{D} \setminus B_k$ is hyperbolically convex, since the intersection of an arbitrary collection of hyperbolically convex sets is again hyperbolically convex. 

For convenience in the proof, we transfer the problem (via a M\"{o}bius isometry) to the upper half-plane $\mathbb{H}$.
Here $\rho_U$ is the complete hyperbolic metric on $U$.

\begin{lemma}\label{L: hyperbolic convexity}
Let $U$ be a simply-connected subdomain of $\mathbb{H}$, with boundary $\partial U$ in  $\mathbb{C}_\infty$. Suppose that $[-1,1] \subset \partial U$ and $\partial U \setminus (-1,1)$ is connected. Define 
\[B \, = \, \{ z \in U \, : \, \omega(z,(-1,1),U) > 1/2 \}\,.\]  Then $B$ is bounded by $[-1,1]$ together with a $\rho_U$-geodesic $\gamma$. Moreover, $\mathbb{H} \setminus B$ is hyperbolically convex in $\mathbb{H}$.
\end{lemma}
\begin{proof}
We may choose a conformal mapping $\varphi: \mathbb{D} \to U$ such that the boundary correspondence maps the semicircle $\{e^{i\theta}\,: \, \theta \in (\pi, 2\pi)\}$ bijectively onto the interval $(-1,1)$. By conformal invariance of harmonic measure, the set $B$ is the image $\varphi(\mathbb{D} \cap \mathbb{H})$ and it is bounded by $[-1,1]$ together with the curve $\varphi((-1,1))$, which is the $\rho_U$-geodesic $\gamma$. 
 Moreover, $\mathbb{H} \setminus \overline{B}$ is connected and simply-connected because $B$ is simply-connected and $B \cup (\mathbb{C}_\infty \setminus \mathbb{H})$ is connected. Write $A = U \setminus \overline{B}$. Let $g$ be a conformal mapping of the domain $\mathbb{H} \setminus \overline{B}$ onto $\mathbb{H}$ such that the analytic boundary arc $\gamma$ of $\mathbb{H} \setminus \overline{B}$ corresponds to the boundary arc $(-1,1)$ of $\mathbb{H}$. Then $g \circ \varphi$ is defined on $\mathbb{D} \cap \mathbb{H}$ and extends continuously to map $(-1,1)$ homeomorphically onto $(-1,1)$. We can therefore define the Schwarz reflection of $g \circ \varphi$ across $(-1,1)$, which gives a conformal map \[ G: \mathbb{D} \,\to\, g(A) \,\cup\, (-1,1) \,\cup\, g(A)^*\,, \] where $g(A)^*$ is the image of $g(A)$ under complex conjugation. We can now extend $g$ to the whole of $\mathbb{H}$ by defining $g(z) = G(\varphi^{-1}(z))$ for all $z \in U$; this agrees with the original definition on the overlap $A$, so defines an analytic function which is a conformal mapping of $\mathbb{H}$ onto $\mathbb{H} \, \cup\,(-1,1) \, \cup \, g(A)^*$. By the result of J{\o}rgensen \cite{J}, the half-plane $\mathbb{H}$ is a hyperbolically convex subset of any domain that contains it, such as $\mathbb{H} \, \cup\,(-1,1) \, \cup \, g(A)^*$. Pulling back by $g$ we find that $U \setminus \overline{B}$ is a hyperbolically convex subset of $\mathbb{H}$, as required.
\end{proof}

We have now completed the proof of Theorem~\ref{T: main}. Next, we show how Lemma~\ref{L: Brown Flinn estimate in the upper half-plane} and Lemma~\ref{L: hyperbolic convexity} suffice to prove the inequalities of Propositions~\ref{P: conformal reflection}, \ref{P: totally real schlicht functions} and \ref{P: level set of harmonic measure}. These three results describe essentially the same situation. If $g$ is the totally real univalent function in Proposition~\ref{P: totally real schlicht functions}, then $(g(z)-a)/(b-a)$ maps $\mathbb{D} \cap \mathbb{H}$ onto a domain $U \subset \mathbb{H}$ as in Proposition~\ref{P: level set of harmonic measure}, and $g(\Gamma)$ is the curve $\gamma$. In the situation of Proposition~\ref{P: level set of harmonic measure}, the points $-1$ and $1$ are certainly contained in the closure of $\gamma$ in $\mathbb{C}$. It follows that $\gamma$ lands at both of these points since $\Lambda_1(\gamma)$ is finite. In fact by Lemma~\ref{L: Brown Flinn estimate in the upper half-plane} we have $\Lambda_1(\gamma) \leqslant \pi$ with equality if and only if $\gamma$ describes a semicircle with diameter $[-1,1]$. This puts us back in the situation of Proposition~\ref{P: conformal reflection}. 

It remains to show for Propositions~\ref{P: totally real schlicht functions} and \ref{P: level set of harmonic measure} that $\gamma$ can only be a semicircle on diameter $[-1,1]$ if $U = \mathbb{H}$. To this end, let $F$ be the unique conformal mapping from $\mathbb{H}$ onto $U$ that fixes $-1, 0$ and $1$ as boundary values. Define \[\Omega \, = \, \mathbb{C} \setminus ((-\infty,-1] \cup [1, \infty))\,.\]
Then $F$ extends by Schwarz reflection to a conformal mapping  \[F: \,\Omega\, \to\, U \cup (-1,1) \cup U^* \,\subseteq\, \Omega\,.\] 

\begin{lemma}\label{L: gamma lies in the unit disc}
If $U \neq \mathbb{H}$ then $\gamma$ is contained in $\mathbb{D} \cap \mathbb{H}$.
\end{lemma}
\begin{proof}
We give two different proofs, both elementary. The first uses the Lindel\"{o}f principle. Suppose $|z| \le 1$, $z \in \mathbb{H}$. Then $\omega(z, [-1,1], \mathbb{H}) \ge 1/2$, so
\[\tfrac{1}{2} \,\le\, \omega(F(z), F([-1,1]),F(\mathbb{H}))\, =\, \omega(F(z), [-1,1], U) \, < \, \omega(F(z), [-1,1],\mathbb{H})\,\]
so $F(z)$ lies in the set \[\left\{w \in \mathbb{H} \,:\, \omega(w,[-1,1],\mathbb{H}) > \tfrac{1}{2} \right\} \, = \, \mathbb{H} \cap \mathbb{D}\, \, ,\] and the lemma follows.

For the second proof, we use the hyperbolic metric $\rho_\Omega$. We observe that
\[\mathbb{D} = \{w \in \Omega \, : \, \rho_\Omega(w, (-1,1)) < \sinh^{-1}(1) \}\,.\]
The Schwarz-Pick lemma says that $F: \Omega \to \Omega$ strictly decreases distances with respect to $\rho_\Omega$, since $F$ is not surjective. By definition, $F$ maps the set $(-1,1)$ into itself, so it also maps the closed $\sinh^{-1}(1)$-neighbourhood of $(-1,1)$ into the open $\sinh^{-1}(1)$-neighbourhood of $(-1,1)$. That is to say, $F\left(\Omega \cap \overline{\mathbb{D}}\right) \subset \mathbb{D}$. Since the relative boundary of $\Omega \cap \overline{\mathbb{D}}$ in $\Omega$ consists of the curve $\gamma$ and its image under complex conjugation, we have $\gamma \subset \mathbb{D} \cap \mathbb{H}$ as required.
\end{proof}

\end{document}